\documentclass[letterpaper, 10 pt, conference]{ieeeconf}  %

\IEEEoverridecommandlockouts                              %
\makeatletter
\let\NAT@parse\undefined
\makeatother

\usepackage{graphicx}
\usepackage{siunitx}

\usepackage{mystyle}
\usepackage{notation}

\renewcommand{\diag}[1]{\mathrm{diag}(#1)}

\title{\LARGE \bf
Parametric Nonconvex Optimization via Convex Surrogates
}

\author{Renzi Wang, Panagiotis Patrinos, and Alberto Bemporad%
\thanks{This work was funded by the European Union (ERC Advanced Research Grant COMPACT, No. 101141351),
KU Leuven internal funding C14/24/103,
FWO project G033822N.
Views and opinions expressed are however those of the authors only and do not necessarily reflect those of the European Union or the European Research Council. Neither the European Union nor the granting authority can be held responsible for them.
R. Wang and A. Bemporad are with the IMT School for Advanced Studies, Lucca, Italy.
        Email: {\tt\small %
        \href{mailto:renzi.wang@imtlucca.it}{<renzi.wang>},%
        \href{mailto:alberto.bemporad@imtlucca.it}{<alberto.bemporad>}%
        \href{mailto:renzi.wang@imtlucca.it,alberto.bemporad@imtlucca.it}{@imtlucca.it}}
P. Patrinos is with Department of Electrical Engineering \textsc{esat-stadius}, KU Leuven, Leuven, Belgium.
Email: {\tt\small \href{mailto:panos.patrinos@kuleuven.be}{panos.patrinos@kuleuven.be}}%
}}

\begin{document}

\maketitle
\thispagestyle{empty}
\pagestyle{empty}

\begin{abstract}

This paper presents a novel learning-based approach to construct a surrogate problem that approximates a given parametric nonconvex optimization problem.
The surrogate function is designed to be the minimum of a finite set of functions,
given by the composition of convex and monotonic terms, so that
the surrogate problem can be solved directly through parallel convex optimization.
As a proof of concept, numerical experiments on a nonconvex path tracking problem confirm the approximation quality of the proposed method.

\end{abstract}

\section{Introduction}
Optimization problems are common in many areas of control systems,
such as robotic motion planning and autonomous driving, as well as in different other engineering fields.
However,
due to the inherent complexity of the problem structure,
the optimization problem is often nonconvex,
posing significant challenges to solve it efficiently.
In practice,
these optimization problems are often solved repeatedly with varying problem parameters,
such as different environment configurations or, in the case of problems arising from optimal control formulations, different initial conditions.
This characteristic gives rise to \emph{parametric optimization problems},
denoting a family of similar optimization problems that differ only in (some of) the problem parameters.

Most of the parametric optimization literature has focused on studying the properties of {\it convex} parametric optimization problems~\cite{MR64,Fia83},
and explicit solutions have been derived for special classes of problems,
such as multiparametric linear~\cite{GN72,BBM03} and (piecewise) quadratic programs~\cite{BMDP02a,TJB03,GBN11,PS10,patrinos2011convex},
to simplify solving different problem instances at run time. However, solving {\it nonconvex} parametric optimization instances repeatedly, 
often under real-time constraints, obtaining fast, accurate, and possibly global solutions remains a crucial challenge in applications, the computational cost per solve often being a critical bottleneck. This opens the door to data-driven approaches, which can exploit the shared structure across problem instances to accelerate solving. 

Data-driven methods to approach parametric optimization problems are broadly referred to as
\emph{learning to optimize} \cite{chen2022learning}, and \emph{amortized optimization} \cite{amos2023tutorial}.
A major category of data-driven approaches is to learn a solution mapping
that directly maps the problem parameters to the optimal solution.
Such a method is often referred to as a {\it fully-amortized} approach \cite[Definition 2]{amos2023tutorial}.
A natural training objective for learning the solution mapping is to minimize the mean squared error between the learned solution and the optimal solution obtained by a solver.
This requires calling a solver offline to generate the training data, which can be computationally expensive, especially for nonconvex optimization problems with several parameters.
In addition, due to the approximation error of the learned model,
the learned solution mapping may not satisfy the constraints of the original problem,
leading to infeasible or suboptimal solutions, unless specific model structures and/or training objectives are designed to ensure, or at least promote, feasibility.

To address the feasibility concern,
two main approaches have been proposed in the literature.
The first approach is to use the learned solution as a warm start for a solver,
which then safeguards feasibility while reducing the number of iterations required \cite{sambharya2024learning}.
However, this approach still requires expensive data collection to train the solution mapping.
The second approach is to learn a solution mapping that explicitly incorporates the constraints by design,
ensuring feasibility of the learned solution \cite{donti2021dc,grontas2026pinet,tordesillas2023rayen}.
However,
these approaches typically train the solution mapping by minimizing the original problem loss over the network parameters via first-order optimization methods,
which effectively reduces to solving the original problem through a reparameterization.
This training procedure may not fully exploit the structure of the original problem or the efficiency of dedicated solvers.

A related but distinct direction places the solver itself at the center,
aiming to learn the algorithmic behavior of the solver rather than the solution mapping.
This line of research is achieved by learning the hyperparameters of the solver,
such as through reinforcement learning \cite{ichnowski2021accelerating}
or by unrolling the algorithm with a feedforward network \cite{oshin2026deep}.
Although most of the works under this category focus on convex optimization problems,
the outcome can be used as a submodule within a nonconvex optimization scheme \cite{oshin2026deep}.
Nevertheless, algorithm unrolling methods still require a solver to generate training data, leaving the computational burden of data collection unresolved.

Instead of learning the solution mapping or the algorithmic behavior of the solver,
another approach is to learn a surrogate problem that approximates the original problem.
A classical way is to construct such an approximation through analytical methods \cite{langeMMOptimizationAlgorithms}.
Recently, data-driven approaches have been proposed to learn a surrogate problem directly from data.
One direction is to learn part of the solution mapping to reduce the original problem
to a simpler one over the remaining variables \cite{bertsimas2022online},
though this still requires accessing the solver to generate training data.
Another direction is to learn a low-dimensional representation of the original problem
through a linear change of variables \cite{NEURIPS2020_6d0c9328}.
However, a linear mapping preserves the nonconvex structure of the problem,
leaving the computational difficulty fundamentally unchanged.
A further direction seeks to learn a coordinate transformation that convexifies the original problem,
such as through a Koopman operator \cite{wu2026koopman,zhang2024latent},
but these methods do not account for the parametric nature of the problem.
The surrogate-based approach also connects to meta-learning,
where the goal is to learn a model that can be quickly adapted to new tasks \cite{finn2017model,rusu2018metalearning}.

\textbf{Contribution}:
The proposed surrogate formulation addresses several limitations identified above.
First, unlike solution mapping and algorithm unrolling approaches,
our method does not require a large corpus of solver-generated (and, possibly, global) solutions as a training prerequisite.
Such data can be optionally incorporated as regularization to improve approximation quality,
but is not a necessity.
Second, by introducing a finite set of convex functions among the components of the surrogate problem,
our approach can exploit the efficiency of dedicated convex solvers,
overcoming the limitation of constrained neural network approaches that rely on general first-order methods.
To summarize, our contribution is twofold:
\begin{inlinelist}
    \item [$i$)] We propose a novel learning-based approach that parametrizes the surrogate problem for parametric nonconvex optimization with quasiconvex components, given by the composition of convex and monotonic functions. %
    The surrogate structure enables us to leverage the powerful tools in convex optimization for solving the surrogate problem in parallel.
    \item [$ii$)] As a proof of concept, numerical experiments on a nonconvex path tracking problem demonstrate the approximation quality of the proposed method.
    The surrogate solution used as an initial guess for original nonconvex problem leads to faster convergence.
\end{inlinelist}

\section{Problem Formulation}
Consider a parametric optimization problem of the form
\begin{equation}\label{eq: optimization_problem}
    \begin{aligned}
    \minimize_{x \in \mathcal{X}}& f(x, p)\\
    \stt\, g(x, p) &\leq 0\\
    \end{aligned}
\end{equation}
where $x\in \re^{\nvar}$ is the decision vector,
and $p \in \re^{\np}$ is the vector of problem parameters.
The function $f:\re^{\nvar} \times \re^{\np} \to \re$ is the nonconvex loss function,
and $g:\re^{\nvar} \times \re^{\np} \to \re^m$ is the constraint function,
with $g(\cdot, p)$ elementwise convex for any fixed $p$.
The domain $\mathcal{X} \subset \re^{\nvar}$ is a convex and compact set.

Since $f$ is nonconvex,
solving the problem \eqref{eq: optimization_problem} for a given parameter $p$ is challenging in general.
We therefore propose to learn a surrogate function $\hat{f}$ that approximates $f$ and admits a structure that allows for efficient optimization.
Specifically, we parameterize $\hat{f}$ in the form
\begin{equation}\label{eq: exact_min}
    \hat{f}(x, p) = \min_{i=1,\ldots,K}\; h_i\big(\bar{f}_i(x, p)\big)
\end{equation}
where ${\bar{f}}_i:\re^{\nvar} \times \re^{\np} \to \re$ is a convex function in $x$ for any fixed $p$ for all $i = 1, \ldots, K$,
and $h_i:\re \to \re$ is a monotonically increasing function for all $i = 1, \ldots, K$.
This parameterization offers two expressiveness properties.
First,
by the composition rule of quasiconvex functions \cite[\S 3.4.4]{boyd2004convex},
each function $x \mapsto h_i(\bar{f}_i(x, p))$ belongs to the strictly broader class of quasiconvex functions in $x$ for any fixed $p$.
Because $h_i$ need not be convex,
$h_i \circ \bar{f}_i$ is not necessarily convex.
Second,
the pointwise minimum of $K$ quasiconvex functions can represent multimodal loss landscapes that no quasiconvex function can capture.

The form \eqref{eq: exact_min} admits a decomposition that
reduces minimizing $\hat{f}$ over the constraints into two stages:
\begin{inlinelist}
    \item [$i$)] Solving $K$ convex optimization problems in parallel by minimizing $\bar{f}_i$ over the constraints, then
    \item [$ii$)] Selecting the optimal solution among the $K$ resulting optimal values of $h_i \circ \bar{f}_i$.
\end{inlinelist}
Note that by minimizing $\bar{f}_i$ we attain the global minimum of each subproblem.
Concretely,
let
\begin{equation}\label{eq: decomposition}
    \begin{aligned}
        \bar{F}_i(p) & = \inf_{x \in \mathcal{X}} \big\{\bar{f}_i (x, p) \mid g(x, p) \leq 0 \big\}
        \quad \text{for } i = 1, \ldots, K, \\
        i^\star(p) &\in \argmin_{i = 1, \ldots, K} h_i\big(\bar{F}_i(p)\big),\\
        x^\star &\in \argmin_{x \in \mathcal{X}} \big\{\bar{f}_{i^\star(p)}(x, p) \mid g(x, p) \leq 0\big\}.
    \end{aligned}
\end{equation}
Then $x^\star \in \argmin_{x\in \mathcal{X}, g(x, p) \leq 0} \hat{f}(x, p)$.

\subsection{Choice of the function}\label{sec: choice_of_function}
Each monotonically increasing function $h_i: \re \to \re$ can be obtained by using a feedforward neural network with $L$ layers,
with nonnegative weights $W_\ell$ and arbitrary biases $b_\ell$,
and monotonic activation functions $\sigma_\ell$ such as ReLU, softplus, and sigmoid for each layer $\ell = 1, \ldots, L - 1$,
and a linear activation function for the last layer $L$.
Alternatively, one can use monotonic rational-quadratic splines proposed in \cite{durkan2019neural}.
Each convex function $\bar{f}_i$ belongs to a chosen parametric family whose coefficients are mappings of $p$.
We collect the parameters of these mappings and the network parameters of $h_i$ into the parameter vector $\Theta$.
Below we list some candidate parametric families for $\bar{f}_i$,
where we omit the index $i$ for notational simplicity.
In all examples,
the coefficient functions are neural networks in $p$,
so $\Theta$ collects their weights and fully determines $\hat{f}_\Theta$.
Any combination of the following families are also a valid choice.

\begin{example}[Parametric quadratic function]%
\[
        \bar{f}(x, p) = \alpha \norm{x}_2^2 + \norm{L(p) x}_2^2 + c(p)^\top x + d(p)
\]
where $\alpha \geq 0$ is a user-specified hyperparameter,
$L(p)$ is a learned lower-triangular matrix,
and $c(p)$, $d(p)$ are learned vector and scalar functions of $p$, respectively.
The lower-triangular structure of $L(p)$ characterizes the positive semidefinite matrix $L^\top(p)L(p)$
with $\nicefrac{\nvar(\nvar+1)}{2}$ learnable parameters,
avoiding the $\nvar^2$ parameters required by a general full matrix
while retaining full rank.
The hyperparameter term $\alpha \norm{x}_2^2$ is kept separate from $L(p)$
to distinguish the fixed regularization $\alpha$ from the learnable component $L(p)$.
Setting $\alpha > 0$ guarantees strong convexity,
ensuring a unique minimizer when the learned problem is solved.
\end{example}
\begin{example}[Max-affine function]%
\[
        \bar{f}(x, p) = \max_{t=1, \dots, L} (A_t(p)x + b_t(p))
\]
where $A_t(p)$, $b_t(p)$ are matrix and vector functions of $p$, respectively.
\end{example}

\begin{example}[Max-squared function]\label{example: max-squared}%
\[
        \bar{f}(x, p) = \sum_{t=1}^L \max(A_t(p) x - b_t(p),\, 0)^2
\]
where $A_t(p)$, $b_t(p)$ are vector and scalar functions of $p$, respectively.
\end{example}

\begin{example}[Input-convex neural network (ICNN) \cite{pmlr-v70-amos17b}]\label{example: icnn}%
To incorporate the parameters, we specifically consider the convex function is of the form
\[
    \bar{f}(x, p) = \ell(x, \varphi(p)),
\]
where $\ell:\re^{\nvar} \times \re^{\nq} \to \re$ is an input-convex neural network that is convex in $x$ for any fixed $\varphi(p)$,
and $\varphi:\re^{\np} \to \re^{\nq}$ is a neural network that maps the problem parameter $p$ to the function parameter $\varphi(p)$,
as proposed in \cite{schallerLearningParametricConvex2025a}.
\end{example}

\subsection{Smoothed minimum}
When the true function $f$ is differentiable,
it would be desirable to ensure that the surrogate function $\hat{f}$ is also differentiable,
so that its gradient can be matched to the true gradient of $f$.
To this end, \emph{during training only}, we consider representing the exact minimum in \eqref{eq: exact_min} 
with a smoothed minimum:
\begin{equation}\label{eq: surrogate_function}
    \hat{f}(x, p) = -\gamma \lse\Big(- \tfrac{1}{\gamma}F\big(x, p\big)\Big)
\end{equation}
where the mapping $F(x, p) = (\bar{f}_1(x, p), \ldots, \bar{f}_K(x, p))$
The function $\lse:\re^K \to \re$ is the log-sum-exponential function defined as 
$\lse(z) = \log\big(\sum_{i=1}^K \exp(z_i)\big)$.
Since $-\lse(-z)$ is a smooth approximation of $\min(z)$,
the surrogate \eqref{eq: surrogate_function} converges to \eqref{eq: exact_min} as $\gamma \to 0$,
and the hyperparameter $\gamma > 0$ trades off approximation quality against numerical conditioning.
After training, i.e., at solve time, one recovers the exact decomposition \eqref{eq: decomposition}
and solves the $K$ convex subproblems in parallel.

\section{Learning the Surrogate Problem}
As mentioned above, many learning-based approaches to parametric nonconvex optimization focus on learning 
a surrogate for the solution $x^\star(p)$,
which needs a solver to produce optimal solutions as supervised training targets.
Such setup introduces a potentially expensive data collection process.
In contrast, in this work we propose to learn the \emph{loss landscape} $f(x, p)$
via a regression problem,
which only requires evaluations of $f$ at arbitrary $(x, p)$ pair, and only optionally
optimal solutions to improve the quality of the surrogate around minima.
Such data is far cheaper to obtain and does not presuppose access to a solver.
We describe the surrogate learning problem in detail below.

Given the dataset $\mathcal{D} =\{x_k, f_k, p_k\}_{k=1}^N$,
where $f_k = f(x_k, p_k)$ denotes the true loss evaluated at the data point $x_k$ with problem parameter $p_k$,
we aim to fit the surrogate $\hat{f}_\Theta$ by minimizing the composite loss:
\begin{equation}\label{eq: total_loss}
    \minimize_{\Theta} \mathcal{L}_{\mathrm{total}}(\Theta)
    \dfn \mathcal{L}(\Theta) + \mathcal{R}(\Theta)
\end{equation}
where each term is described next.

\subsection{Data fitting loss}
The primary objective is to minimize the mean squared error between the true loss and the approximated loss:
\begin{equation}\label{eq: loss_function}
    \mathcal{L}(\Theta) = \frac{1}{N} \sum_{k=1}^N \big(f_k - \hat{f}_\Theta(x_k, p_k)\big)^2.
\end{equation}

\subsection{Regularization}
\paragraph{Optimality regularization (optional)}
Let $\mathcal{D}^\star = \{x_i^\star, f_i^\star, p_i^\star, \lambda_i^\star\}_{i=1}^{M_1} \subseteq \mathcal{D}$
be the subset of the training data at which $x_i^\star$ is the optimal solution of the optimization problem \eqref{eq: optimization_problem} with parameter $p_i^\star$, where we additionally assume access to the dual variables $\lambda_i^\star \in \re^m$ associated with the constraints $g(x, p) \leq 0$. This are often immediately available when the data is collected from a numerical solver.  

At the optimal solution $x_i^\star$,
the stationarity condition of the Lagrangian function requires
\begin{equation}
    \nabla_x f(x_i^\star, p_i^\star) + \nabla_x g(x_i^\star, p_i^\star)\lambda_i^\star = 0.
\end{equation}
To encourage the surrogate function to respect the optimality condition at these points,
we consider the following regularization term:
\begin{equation}\label{eq: regularization_optimality}
    \mathcal{R}_1(\Theta) = \frac{w_1}{M_1} \sum_{i=1}^{M_1} \norm{\nabla_x \hat{f}_\Theta(x_i^\star, p_i^\star) +  \nabla_x g(x_i^\star, p_i^\star)\lambda_i^\star }^2_2
\end{equation}
where $w_1 > 0$ is a hyperparameter that controls the weight of such a regularization term.

Importantly, we remark that
points where the solver converged only within a loose numerical tolerance
are still valuable information for training the surrogate function and 
can still be included in $\mathcal{D} \setminus \mathcal{D}^\star$.

\paragraph{Gradient matching regularization}
Let $\mathcal{D}^{\dagger} = \{x_i, f_i, p_i\}_{i=1}^{M_2} \subseteq \mathcal{D}$
be the subset of the training data at which the true gradient $\nabla_x f_i \dfn \nabla_x f(x_i, p_i)$ is directly accessible.
In this case, one can match the surrogate gradient to the true gradient by considering the following regularization term:
\begin{equation}\label{eq: regularization_curvature}
    \mathcal{R}_2(\Theta) = \frac{w_2}{M_2} \sum_{i=1}^{M_2} \norm{\nabla_x \hat{f}_\Theta(x_i, p_i) - \nabla_x f_i}^2_2
\end{equation}
where $w_2 > 0$ is a hyperparameter that controls the weight of the regularization term.
Unlike \eqref{eq: regularization_optimality},
this regularization is applicable at \emph{any} point where the gradient is available.
Clearly, if $\mathcal{D}^\star = \mathcal{D}^{\dagger}$ and the Karush-Kuhn-Tucker (KKT) optimality conditions hold {\it exactly} at the optimal points, i.e.,
\[
    \nabla_x f_i^\star + \nabla_x g(x_i^\star, p_i^\star)\lambda_i^\star = 0
\]
the regularization term \eqref{eq: regularization_optimality} is equivalent to the regularization term \eqref{eq: regularization_curvature}.

\subsection{Projected sampling}\label{sec: projected_sampling}
A natural approach to construct the dataset $\mathcal{D}$ is to sample the data points $\{x_k\}_{k=1}^N$ uniformly from the domain $\mathcal{X}$.
However, uniform sampling tends to under-represent the boundary of $\mathcal{X}$,
which is critical for capturing the behavior of the loss function near the boundary and ensuring the surrogate function is accurate in these regions.
To ensure adequate coverage of the boundary of the domain,
we employ a projected sampling strategy:
\begin{enumerate}[label=\roman*)]
    \item [$i$)] Construct an enlarged domain $\tilde{\mathcal{X}}$ such that $\mathcal{X} \subseteq \tilde{\mathcal{X}}$.
    \item [$ii$)] Sample the data points $\{x_k\}_{k=1}^N$ uniformly from the enlarged domain $\tilde{\mathcal{X}}$.
    \item [$iii$)] Project the sampled data points onto the original domain $\mathcal{X}$ by solving a projection problem:
    \begin{equation}
        x_k^{\mathrm{proj}} = \arg\min_{x \in \mathcal{X}} \norm{x - x_k}^2_2.
    \end{equation}
\end{enumerate}
The projection is the identity for points already inside $\mathcal{X}$ and maps boundary-violating samples to the nearest feasible point,
thereby enriching the dataset near the boundary.
Notably, the projection enforces only $\mathcal{X}$
and deliberately does not enforce the constraints $g(x, p) \leq 0$.
This is because $f$ is well-defined over the entire domain $\mathcal{X}$,
and the feasible set $\mathcal{X}_{p_i} \dfn \{x\in\mathcal{X} \mid g(x, p_i) \leq 0\} \subseteq \mathcal{X}$ varies with the parameter $p$.
By sampling from the larger set $\mathcal{X}$,
the dataset provides the coverage of the full operational feasible region across all parameter values,
yielding a surrogate capturing the global structure of the loss landscape.
Samples restricted to $\mathcal{X}_{p_i}$ only provide local and parameter-specific information,
making the surrogate difficult to generalize well beyond $\mathcal{X}_{p_i}$.
Empirically, restricting samples to $\mathcal{X}_{p_i}$ degrades the performance,
which is consistent with the intuition above.

For parameters $p_k$,
two approaches can be considered.
The first approach is to sample $p_k$ uniformly from a predefined parameter set $\mathcal{P} \subseteq \re^{\np}$.
The second is applicable when prior solution data is readily available,
such as from a model predictive controller that has been deployed on a system,
producing a dataset $(x^\star_i, p_i)_{i=1}^M$.
In this case,
the parameter values can be directly inherited from this dataset.
The projected sampling strategy is then applied to augment the coverage of the dataset in the $x$-space,
i.e., for each fixed parameter value $p_i$,
additional samples $(x_j^{\mathrm{proj}})$ are generated and paired with $p_i$ to enrich the dataset.
The true loss $f_k = f(x_k^{\mathrm{proj}}, p_k)$ is then evaluated at the projected points,
which form the final dataset $\mathcal{D}$.

\section{Numerical Experiments}
\subsection{Static function approximation}
\begin{figure}[!tb]
    \centering
    \begin{subfigure}{0.23\textwidth}
        \centering
        \includegraphics[width=\textwidth]{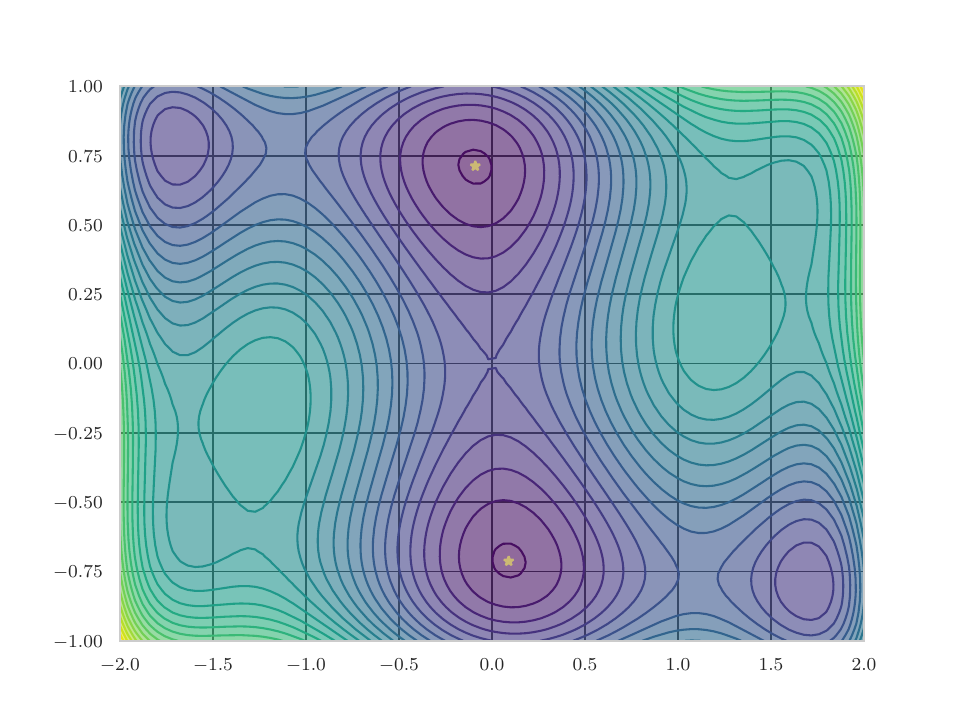}%
        \caption{Ground truth}
    \end{subfigure}
    \hfill
    \begin{subfigure}{0.23\textwidth}
        \centering
        \includegraphics[width=\textwidth]{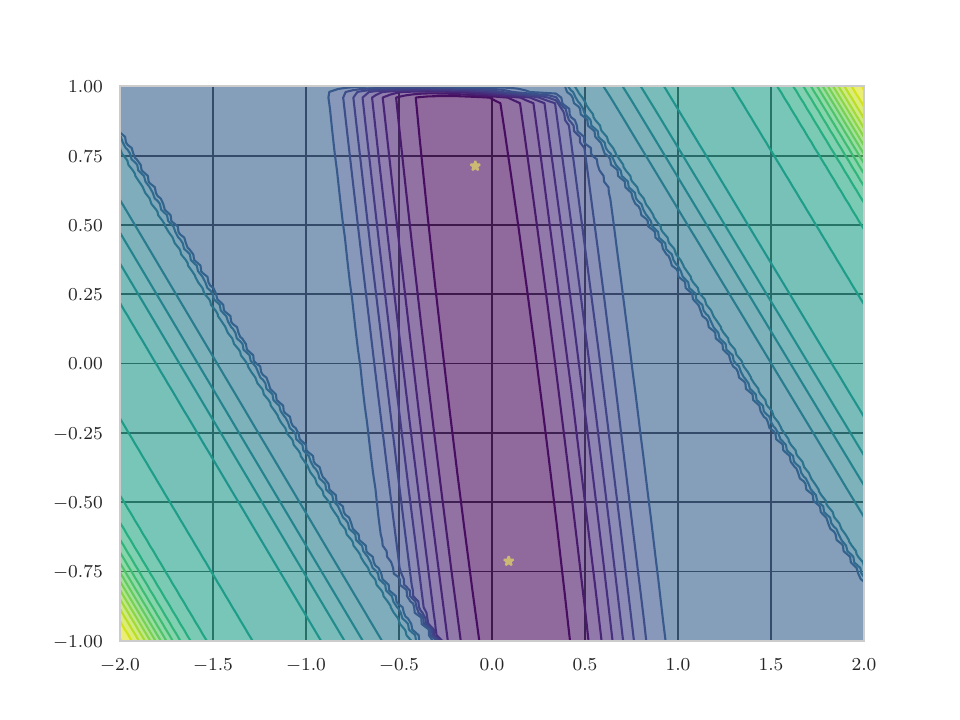}%
        \caption{$K = 1$}
    \end{subfigure}
    \vfill
    \begin{subfigure}{0.23\textwidth}
        \centering
        \includegraphics[width=\textwidth]{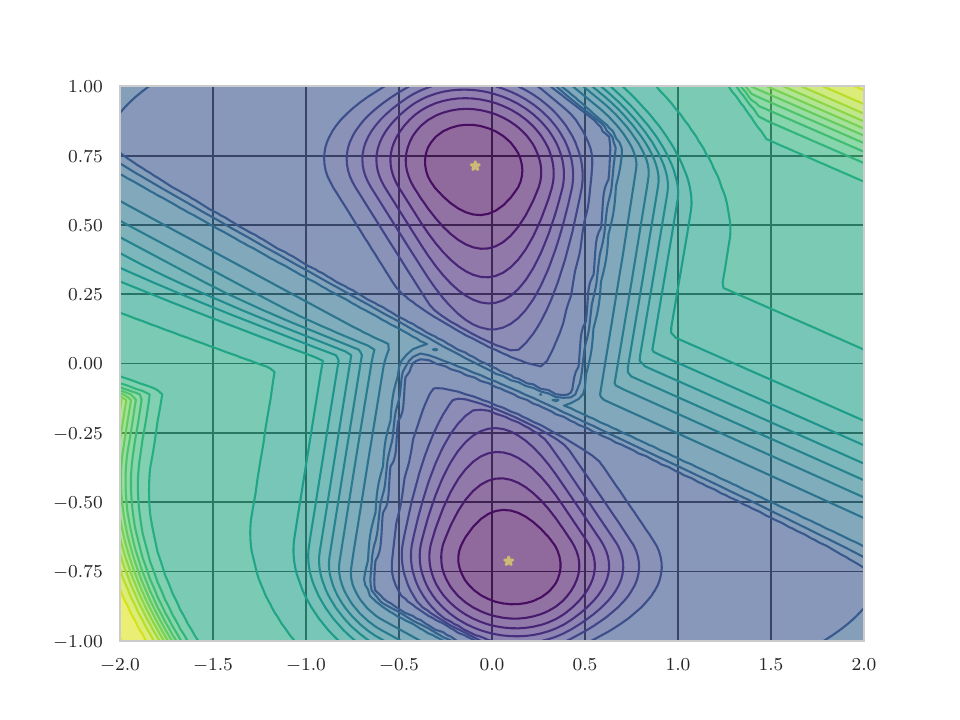}%
        \caption{$K = 2$}
    \end{subfigure}
    \hfill
    \begin{subfigure}{0.23\textwidth}
        \centering
        \includegraphics[width=\textwidth]{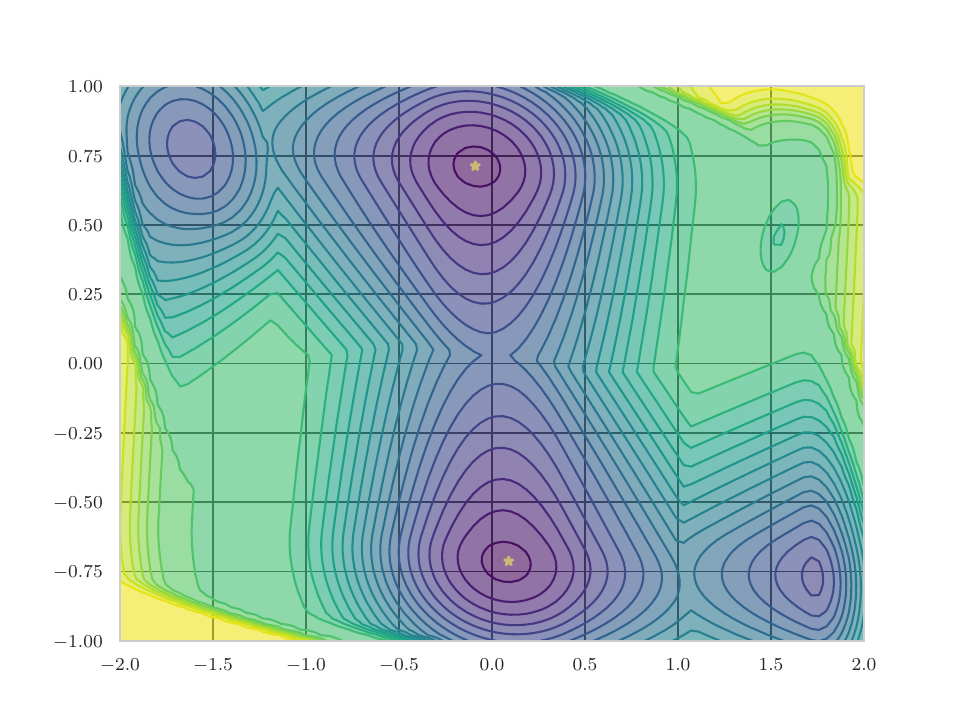}%
        \caption{$K = 5$}
    \end{subfigure}
    \caption{Level sets of the ground truth function and the surrogate functions with different numbers of components $K$,
    where the dark color indicates the low function value and the yellow stars indicate the global minima.}
    \label{fig:camel-six-hump}
\end{figure}

We use the six-hump camel back function \cite{molga2005test} as a test function
to illustrate that a nonconvex function can be effectively approximated by the proposed formulation \eqref{eq: exact_min}.
The function is defined as
\[
    f(x_1, x_2) = (4 - 2.1x_1^2 + \frac{x_1^4}{3})x_1^2 + x_1x_2 + (4x_2^2 - 4)x_2^2.
\]
This function has two global minima equal to $f(x) = -1.0316$,
located at $(x_1, x_2) = \bmat{0.0898 & -0.7126}$
and $(x_1, x_2) = \bmat{-0.0898 & 0.7126}$.
We sample $N=1000$ samples uniformly in the domain $x_1 \in [-2, 2]$ and $x_2 \in [-1, 1]$.
The monotonic function $h_i$ is parameterized as a two-layer neural network with 5 and 3 neurons, respectively.
Each layer is activated by the $\mathtt{tanh}$ function.
To simplify the training,
we enforce that $h_i$ are shared across the $K$ components of the surrogate function.
The convex functions are parameterized with the max-squared form presented in \cref{example: max-squared}.
Particularly,
we set $L = 10$.
As this function has no parameter $p$,
the matrices $A_t$ and the vectors $b_t$ are directly learned as free parameters, without any dependence on $p$.
The surrogate function is implemented through \texttt{jax-sysid} \cite{Bem25}.
The model is trained for 1000 epochs using \texttt{Adam} \cite{kingma2014adam} with a learning rate of $10^{-3}$,
followed by 5000 epochs of fine-tuning with \texttt{L-BFGS} \cite{byrd1995limited}.
\cref{fig:camel-six-hump} illustrates the level sets of the ground truth function and the surrogate functions with different component number $K$.

As demonstrated in \cref{fig:camel-six-hump},
increasing the number of components $K$ improves approximation quality.
With $K=2$ components,
the surrogate function already captures the region with two global minima.
By further increasing $K$ to 5,
the surrogate function more closely approximates the overall landscape of the true function.

\begin{figure}[!tb]
    \centering
    \includegraphics[width=0.45\textwidth]{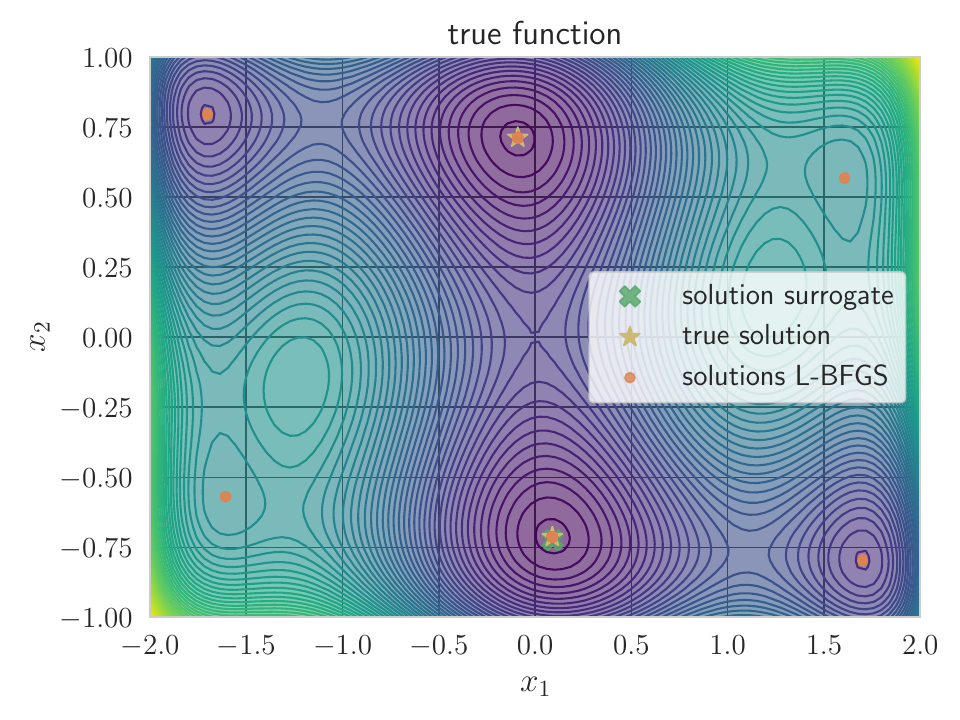}%
    \caption{Solutions obtained by minimizing the surrogate directly
    and by running \texttt{L-BFGS} from 100 random initializations on the original function}
    \label{fig:camel-six-hump-convergence}
\end{figure}

To further evaluate the quality of the surrogate function,
we optimize the surrogate function using the procedure described in \eqref{eq: decomposition} with $K=5$.
As a baseline, we also solve the original optimization problem using \texttt{L-BFGS} with 100 random initializations.
In particular, 64 out of 100 random initializations of \texttt{L-BFGS} converge to the global minima,
while the remaining 36 converge to local stationary points,
as illustrated in \cref{fig:camel-six-hump-convergence}.
In contrast, the surrogate solution $x^\star = \bmat{0.092 & -0.724}$ is close to a global minimizer.
When \texttt{L-BFGS} is initialized at $x^\star$,
it converged to the global minimum at $\bmat{0.0898 &  -0.7126}$ with 5 iterations.
The result demonstrates that learning a surrogate approximating the original function
effectively captures the global landscape of the original function,
such that optimizing the surrogate function directly yields a solution close to a global minimum.

\subsection{Optimal control problem approximation}
We consider the following simple unicycle model for a mobile robot
\begin{equation}
\bmat{\dot{p}_x \\ \dot{p}_y \\ \dot{\psi}} = \bmat{v\cos{\psi} \\ v\sin{\psi} \\ \omega}
\end{equation}
where $(p_x, p_y) \in \re^2$ is the Cartesian coordinate of the robot and $\psi \in (-\pi, \pi]$  is the orientation in the global frame.
The control input $u = (v, \omega) \in \re^2$ is the longitudinal velocity and the angular velocity of the robot, respectively,
with physical constraints $v \in [-v_{\max}, v_{\max}]$ and $\omega \in [-\omega_{\max}, \omega_{\max}]$.
Specifically, we let $v_{\max} = \SI{1.2}{\meter\per\second}$ and $\omega_{\max} = \nicefrac{\pi}{3}\,\SI{}{\radian\per\second}$.
The objective is to track a Lissajous curve while maintaining a lateral deviation within half the track width, i.e., $d \leq \bar{d}/2 = \SI{0.3}{\meter}$,
where $d$ is the lateral deviation from the track centerline, defined in the Frenet frame presented below.
The Lissajous curve is defined as:
\begin{align*}
p_x =& A\sin(at + \delta) \\
p_y =& B\sin(bt)
\end{align*}
with $A = 1.5$, $B = 2$, $a = 3$, $b = 2$, and $\delta = \frac{\pi}{2}$.
We generate the reference trajectory by sampling the Lissajous curve with a reference velocity of $\SI{0.72}{\meter\per\second}$, and the sampling time is $\dd t = \SI{0.1}{\second}$.

\begin{figure}[!bt]
    \centering
    \includegraphics[width=0.35\textwidth]{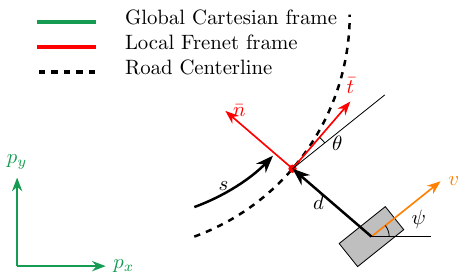}
    \caption{Illustration of the Frenet frame.}
    \label{fig: frenet_coordinates}
\end{figure}
To reduce the dimension of the problem parameter,
we convert the system state from the Cartesian coordinate to the Frenet frame \cite{qian2016hierarchical}.
Frenet coordinates describe the robot pose w.r.t. the centerline of the road,
which is illustrated in \cref{fig: frenet_coordinates}.
The Frenet coordinate system consists of the arc length $s$,
which represents the travel distance along the road;
the lateral offset $d$,
and the heading error $\theta$ between the vehicle yaw angle and the heading of the road.
We denote the system state in the Frenet frame as $x = (s, d, \theta) \in \re^3$,
where we normalize the longitudinal displacement so that $s \in [0, 1]$.
Because the unicycle model is input-affine,
we write the discrete time dynamics as
\begin{equation}
x_{k+1} = f(x_k, u_k) \dfn A(x_k) + B(x_k) u_k
\end{equation}
where $A: \re^3 \to \re^3$ and $B: \re^3 \to \re^{3\times 2}$
are the state-dependent drift and control matrices, respectively.
Since the reference trajectory is the center line of the track with forward velocity $\SI{0.72}{\meter\per\second}$,
the reference longitudinal displacement and the lateral displacement are always fixed among different problems.
The only difference is the curvature of the reference trajectory at different time steps,
which is determined by the curvature of the Lissajous curve.

We apply a model predictive control (MPC) approach to solve the tracking problem.
We set the prediction horizon $N_{p} = 10$
and consider a move blocking strategy with the move blocking horizon $N_b = 4$.
The optimal control problem is defined as:
\begin{subequations}\label{eq: ocp}
    \begin{align}
        \min_{u_0, \ldots, u_{N_p-1}} &
        \tlsum_{k=0}^{N_p-1} \!\!
        \norm{x_k {-} x_k^{\mathrm{ref}}}_{Q}^2
        {+} \norm{u_k}_{R}^2 {+} \norm{x_{N_p} {-} x_{N_p}^{\mathrm{ref}}}_{Q_N}^2 \\
        \textrm{s.t. } & x_0 = x_t \\
        &x_{k+1} = f(x_k, u_k), \quad\forall k \in [0, N_p{-}1] \\
        &u_k\in \mathbb{U}, \quad\forall k \in [0, N_p{-}1] \\
        &u_k = u_{N_b} \quad \forall k \geq N_b, \\
        &A(x_0) + B(x_0) u_0 \in \mathbb{X},
    \end{align}
\end{subequations}
where
$\mathbb{U} \dfn [0, v_{\max}] \times [-\omega_{\max}, \omega_{\max}]$,
and $\mathbb{X} \dfn \{x \in \re^3\mid d \in [-\frac{\bar{d}}{2}, \frac{\bar{d}}{2}]\}$.
Here the velocity is constrained to be non-negative since the robot is expected to only move forward along the track,
and the track boundary constraint is only applied to the first time step to simplify the problem.
The weights of the cost function are set as
$Q = \diag{5L, 3, 0.1}$,
$R = 0.01I$,
$Q_N = \diag{100L, 15, 0.5}$,
with $L = 25.68$ the total length of the track.
The system state at time step $t$ is denoted as $x_t$,
and $x_k^{\mathrm{ref}}$ is the reference state at time step $k$.
As noted above,
the reference in the local Frenet frame is fixed among all problems.
And we give $p=(x_t, (\Delta\psi_k^{\mathrm{ref}})_{k=0}^{N_p})$ for our learning problem,
where $\Delta \psi_k^{\mathrm{ref}} = \psi_k^{\mathrm{ref}} - \psi_0^{\mathrm{ref}}$
is the reference centerline heading angle difference at time step $k$.

\subsubsection{Training data collection}
We sample 1000 problems from the Lissajous path tracking problem \eqref{eq: ocp}
over 100 simulated laps, with 10 problems sampled per lap.
In each lap, the robot is randomly initialized on the Lissajous track
with its lateral displacement uniformly sampled from $[-\bar{d}/2 + 0.1, \bar{d}/2 - 0.1]$
and its heading angle uniformly sampled from $[-\pi/6, \pi/6]$.
To ensure diversity across the 1000 sampled problems,
a small perturbation is applied to the robot state at each time step during data collection,
with lateral displacement perturbed by $d \in [-0.02, 0.02]\,\SI{}{\meter}$
and heading angle perturbed by $\theta \in [-\pi/36, \pi/36]\,\SI{}{\radian}$.
For problem $k$ with parameter $p_k = (x_t, (\Delta\psi_i^{\mathrm{ref}})_{i=0}^{N_p})$,
\eqref{eq: ocp} is solved to collect the optimal solution $\ve{u} = (u_0, \ldots, u_{N_p-1})$.
The gradient of the loss function at this point is also recorded to apply the gradient-matching regularization \eqref{eq: regularization_curvature} during training.
We denote by $\underline{u} = (-v_{\max}, -\omega_{\max})$ and $\overline{u} = (v_{\max}, \omega_{\max})$
the lower and upper bounds of the control input, respectively.
We apply the projected sampling strategy described in \cref{sec: projected_sampling}
with the enlarged domain $[\underline{u} - \Delta u, \overline{u} + \Delta u]^{N_b}$,
where $\Delta u = 0.5 (\overline{u} - \underline{u})$,
to sample additional data points for training the surrogate function.
For each problem,
we generate 300 additional data points.

We instantiate the surrogate \eqref{eq: surrogate_function} with $K = 2$ convex functions,
each parameterized by an ICNN presented in \cref{example: icnn} with 2 hidden layers of 5 units each.
The hidden layers apply the \texttt{softplus} function as the activation function and the output layer is a linear layer.
The \texttt{softplus} activation is preferred over \texttt{ReLU} owing to its smoothness,
which is required for the gradient-matching regularization \eqref{eq: regularization_curvature}.
To ensure nonnegativity of the ICNN weights,
a \texttt{softplus} function is applied to the latent weight variables.
The parameter encoder $\varphi : \re^{\np} \to \re^{\nq}$ is implemented as a separate multi-layer perceptron (MLP) for each layer of the ICNN with 2 hidden layers.
The layer width is dynamically determined by $\nicefrac{\np + \nq}{2}$.
\texttt{softplus} activations on the hidden layers,
and a $\tanh$ output activation.
The monotonic function $h_i$ is simply an identity function in this experiment.
The full network leads to a learnable parameter $\Theta \in \re^{28128}$.

The training is implemented through \texttt{jax-sysid} \cite{Bem25}.
The parameters are optimized by minimizing the composite loss \eqref{eq: total_loss}
for 1000 epochs with the \texttt{Adam} \cite{kingma2014adam},
followed by 2000 iterations with the \texttt{L-BFGS} \cite{byrd1995limited}.
Training is repeated with 4 independent random initializations,
and the model achieving the highest $R^2$ score on the training set is retained.

\subsubsection{Evaluation}
We evaluate the proposed method as a primal initialization strategy for solving problem \eqref{eq: ocp},
where the surrogate solution is used as the initial guess of the primal variable
for \texttt{sqpmethod} in \texttt{CasADi} \cite{Andersson2018} with the quadratic programming (QP) solver \texttt{qpOASES} \cite{Ferreau2014}.
The sequential quadratic programming (SQP) solver then refines the surrogate solution to produce the final control input.
The convex problems in the surrogate problem implemented in \texttt{CVXPY} \cite{agrawal2018rewriting} are solved with \texttt{Clarabel} \cite{Clarabel_2024}.
We run the closed-loop tracking experiment for 400 time steps,
and compare the proposed method against two baselines:
\begin{enumerate}[label=\roman*)]
    \item \texttt{cold-start}: the default zero initialization of \texttt{CasADi};
    \item \texttt{shifted solution}: the solution from the previous time step shifted by one step,
    with the last control input repeated for the final step.
    The initialization is combined to solve the original problem to full convergence and with 2 iterations of SQP,
    denoted as \texttt{shifted solution-full} and \texttt{shifted solution-2}, respectively.
\end{enumerate}
For the proposed method, we report two variants:
\texttt{learning-based-full}, which solves the SQP to full convergence,
and \texttt{learning-based-2}, which limits the SQP to 2 iterations.
To evaluate generalization to unseen initial conditions,
the robot is subjected to a larger perturbation during the closed-loop experiment,
with lateral displacement perturbed by $d \in [-0.05, 0.05]\,\SI{}{\meter}$
and heading angle perturbed by $\theta \in [-\pi/24, \pi/24]\,\SI{}{\radian}$,
ensuring that the states encountered at test time lie outside the training distribution.

\begin{figure}
    \centering
    \includegraphics[width=0.45\textwidth]{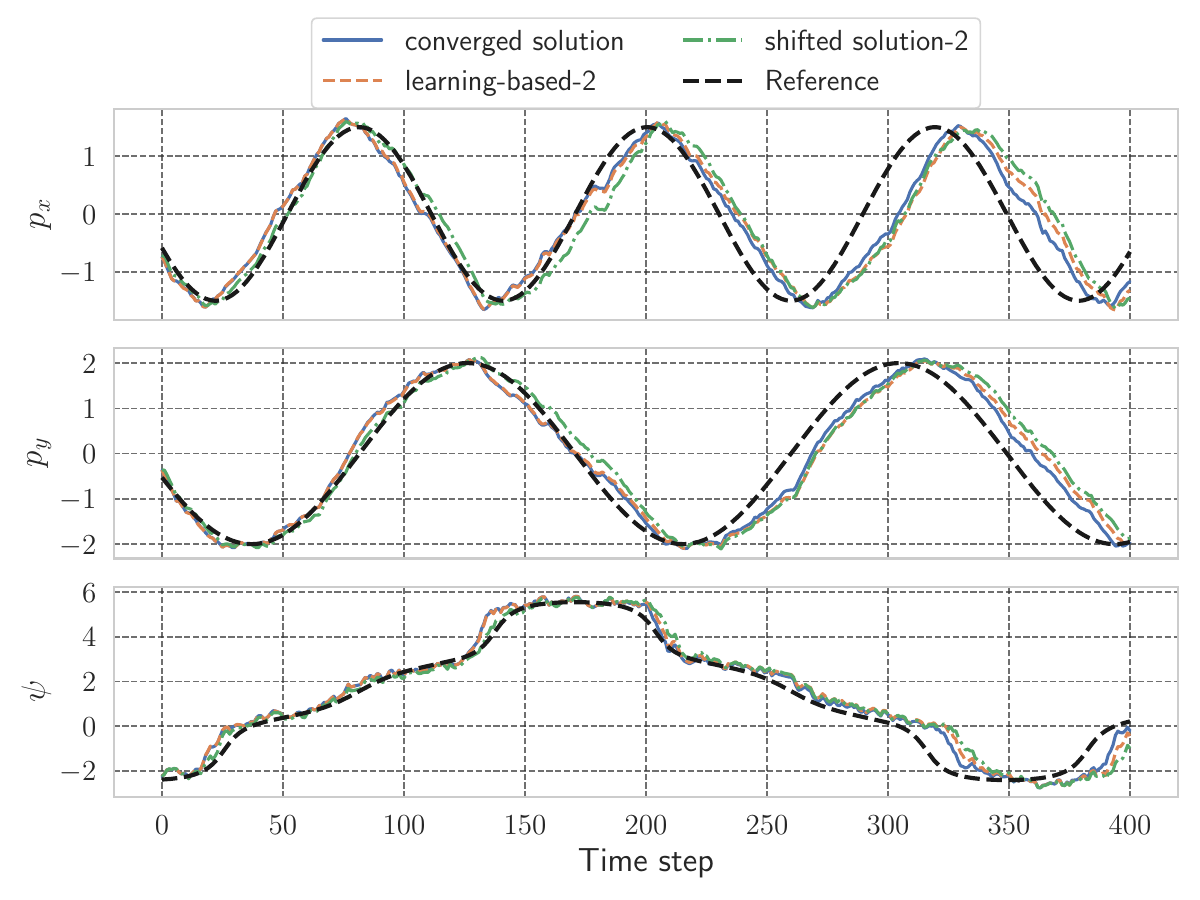}
    \caption{Closed-loop trajectory with different initialization methods.
    The trajectories of the methods \texttt{cold-start}, \texttt{shifted solution-full}, and \texttt{learning-based-full}
    are almost indistinguishable,
    and are denoted as \texttt{converged solution} for clear visualization.
    }
    \label{fig: state_trajectories}
\end{figure}
\cref{fig: state_trajectories} shows the closed-loop trajectory of the robot with different initialization methods.
As all \texttt{cold-start}, \texttt{shifted solution-full}, and \texttt{learning-based-full} methods solve the original problem to full convergence,
leading to negligible differences in the closed-loop trajectory,
we denote the trajectory as \texttt{converged solution} in \cref{fig: state_trajectories} for clear visualization.
From \cref{fig: state_trajectories},
we see that when the solver runs to full convergence,
all methods closely follow the reference trajectory.
Even with only 2 iterations of SQP,
the learning-based method can still achieve a closed-loop trajectory that is close to the full convergence case.
Unlike the trajectory tracking,
where deviations from a time-indexed reference generate a correction signal,
the path tracking formulation re-generates the reference at each time step based on the current state.
So the state displacement induced by the suboptimality does not generate a correction signal,
causing gradual drift.
The results show that 2 SQP iterations are sufficient to keep the state displacement small enough so that the drift remains negligible over the full trajectory.
Compared to the shifted solution,
the learning-based method can achieve a closer trajectory to the reference.
This demonstrates the quality of the surrogate solution as an initialization for the original problem.

\begin{table}
\centering
\begin{tabular}{lrrrr}
\toprule
Method &  mean & min & max & std. \\
\midrule
\texttt{cold-start}              & 4.140 & 1.857 & 15.565 & 1.063 \\
\texttt{shifted solution-full}   & 2.582 & 1.467 & 8.686  & 0.567 \\
\texttt{shifted solution-2}      & 1.717 & 1.151 & 4.426  & 0.382 \\
\texttt{learning-based-full} (total) & 2.739 & 1.654 & 5.472  & 0.329 \\
\texttt{learning-based-2} (total) & 2.211 & 1.691 & 3.089  & 0.171 \\
\midrule
\texttt{learning-based-full} (init.) & 0.618 & 0.504 & 1.230  & 0.087 \\
\texttt{learning-based-2} (init.)    & 0.622 & 0.500 & 1.451  & 0.086 \\
\texttt{learning-based-full} (SQP)    & 2.122 & 1.117 & 4.373  & 0.310 \\
\texttt{learning-based-2} (SQP)      & 1.589 & 1.105 & 2.521  & 0.154 \\
\bottomrule
\end{tabular}
\caption{Computation time (\SI{}{\milli\second}) of different methods.}
\label{tab: computation_time}
\end{table}
\cref{tab: computation_time} reports the per-step computation time for each method.
The second block of \cref{tab: computation_time} breaks down the total runtime of the proposed method
into the initialization time that solves the surrogate problem,
and the SQP solve time.
The runtime for solving the surrogate problem is recorded as the maximum runtime across $K$ convex subproblems.
Compared to the \texttt{cold-start} method,
the learning-based method reduces the total runtime by around $33\%$,
with the runtime for SQP reduced by around $48\%$.
However, the overhead of solving the surrogate problem narrows the advantage
over \texttt{shifted-solution},
which serves as a strong baseline in this experiment for two reasons:
\begin{inlinelist}
    \item [$i$)] the moderate perturbation keeps the system close to deterministic;
    \item [$ii$)] the high-curvature segments of the Lissajous trajectory induce a near bang-bang policy in $\omega$,
where the optimal input remains saturated across consecutive time steps and is largely unchanged by the time shift.
\end{inlinelist}
Noticeably, the worst-case runtime of the learning-based method is significantly lower than the \texttt{shifted solution} method,
demonstrating that the surrogate method is more robust in providing good initializations for the original problem.

\begin{figure}
    \centering
    \begin{subfigure}{0.23\textwidth}
        \centering
        \includegraphics[width=\textwidth]{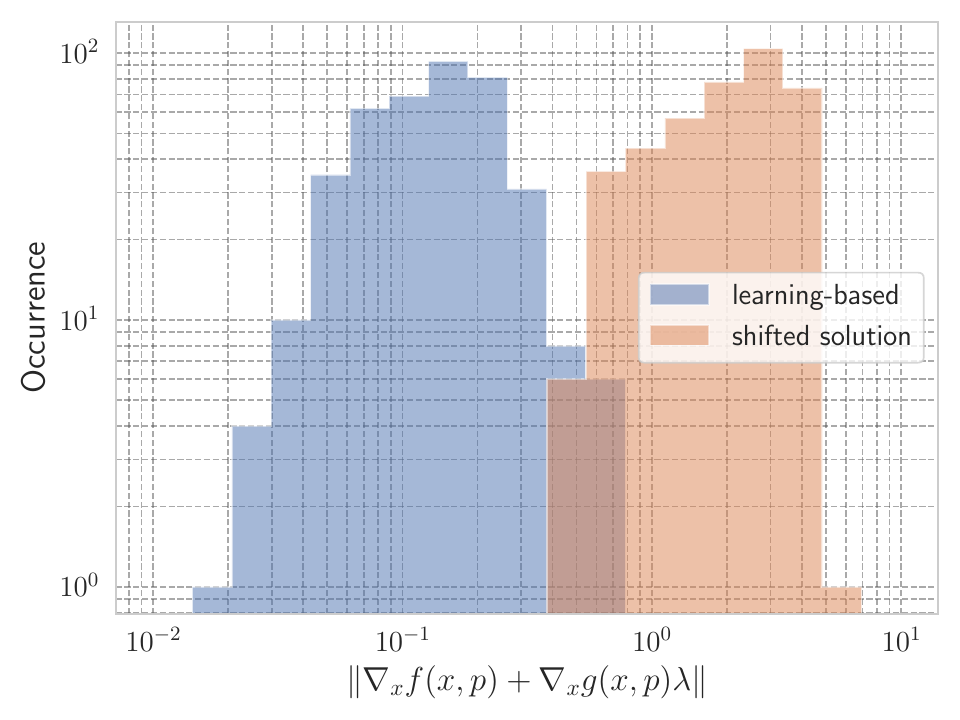}
        \caption{Initial guess}
    \end{subfigure}
    \hfill
    \begin{subfigure}{0.23\textwidth}
        \centering
        \includegraphics[width=\textwidth]{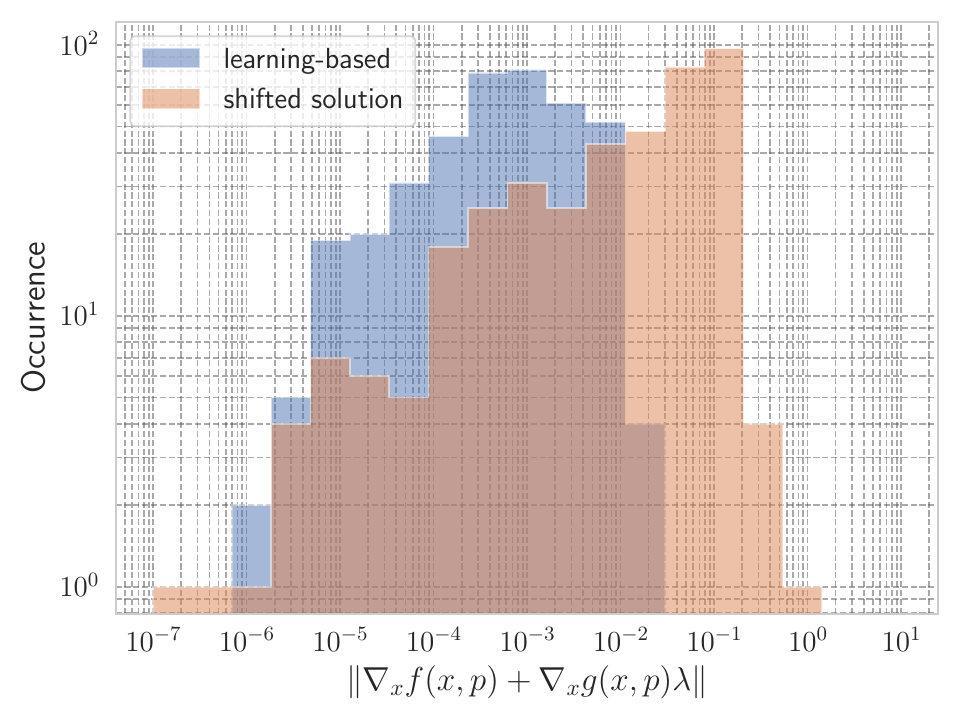}
        \caption{After 2 SQP iterations}
    \end{subfigure}
    \caption{Histogram of the stationarity residual of the original problem at the initial guesses
    and after 2 iterations of SQP, for the learning-based method and the shifted solution.}
    \label{fig: kkt_residual_histogram}
\end{figure}

To further compare the quality of the surrogate solution against the shifted solution,
we examine the stationarity residual
\(
    \norm{\nabla_x f(x, p) + \nabla_x g(x, p) \lambda}
\)
as a solution quality metric,
where $\lambda \geq 0$ is the dual variable for the inequality constraint $g(x, p) \leq 0$.
Since the surrogate solution is primal feasible by construction,
and the equality constraints are eliminated by single shooting,
the stationarity condition is the only remaining KKT condition to satisfy.
For the shifted solution that could be potentially infeasible due to the time shift,
a comprehensive KKT residual would include an additional nonnegative feasibility term,
so the comparison on stationarity residual alone is conservative and understates the advantage of the proposed method.
The evaluation is performed at the initial guess, and after 2 iterations of SQP, respectively.
To obtain a consistent dual variable at the initial guess,
we pass the primal initial guess to \texttt{CasADi} with the iteration limit set to zero,
so that \texttt{CasADi} initializes the dual variable without solving any QP.
After 2 iterations of SQP,
the same residual is evaluated at the updated primal-dual pair.
The histogram of the residual is shown in \cref{fig: kkt_residual_histogram}.
From \cref{fig: kkt_residual_histogram},
we observe that the surrogate solution has a significantly lower stationarity residual value than the shifted solution,
confirming that the surrogate solution is closer to a stationarity point of the original problem.
After 2 iterations of SQP,
the residuals of both approaches reduce,
with the learning-based method retaining a lower residual than the shifted solution.

\section{Conclusion}
In this paper we proposed a novel learning-based approach for constructing surrogate problems of a rather broad class of parametric nonconvex optimization problems.
Crucially, by using the minimum of quasiconvex components to parameterize its objective function, the surrogate problem can be solved directly via parallel convex optimization.
We validated the approximation quality of the proposed formulation in numerical experiments, including a nonconvex path-tracking problem. Our approach can be used either to replace the original problem with the surrogate problem, or to use the latter to provide a good initial guess for solving the original problem, reducing computation time and increasing the likelihood of obtaining global solutions.
Future work includes extending the proposed approach to nonconvex inequality constraints,
incorporating the method into active learning schemes for high-dimensional problems,
embedding the surrogate problem into adaptive control schemes,
and validating the framework on more complex problems.

\addcontentsline{toc}{section}{References}
\bibliographystyle{IEEEtran}
\bibliography{IEEEabrv,reference}

\begin{thebibliography}{10}
\providecommand{\url}[1]{#1}
\csname url@samestyle\endcsname
\providecommand{\newblock}{\relax}
\providecommand{\bibinfo}[2]{#2}
\providecommand{\BIBentrySTDinterwordspacing}{\spaceskip=0pt\relax}
\providecommand{\BIBentryALTinterwordstretchfactor}{4}
\providecommand{\BIBentryALTinterwordspacing}{\spaceskip=\fontdimen2\font plus
\BIBentryALTinterwordstretchfactor\fontdimen3\font minus
  \fontdimen4\font\relax}
\providecommand{\BIBforeignlanguage}[2]{{%
\expandafter\ifx\csname l@#1\endcsname\relax
\typeout{** WARNING: IEEEtran.bst: No hyphenation pattern has been}%
\typeout{** loaded for the language `#1'. Using the pattern for}%
\typeout{** the default language instead.}%
\else
\language=\csname l@#1\endcsname
\fi
#2}}
\providecommand{\BIBdecl}{\relax}
\BIBdecl

\bibitem{MR64}
O.~Mangasarian and J.~Rosen, ``Inequalities for stochastic nonlinear
  programming problems,'' \emph{Operations Research}, vol.~12, pp. 143--154,
  1964.

\bibitem{Fia83}
A.~Fiacco, \emph{Introduction to Sensitivity and Stability Analysis in
  Nonlinear Programming}.\hskip 1em plus 0.5em minus 0.4em\relax London, U.K.:
  Academic Press, 1983.

\bibitem{GN72}
T.~Gal and J.~Nedoma, ``Multiparametric linear programming,'' \emph{Management
  Science}, vol.~18, no.~7, pp. 406--422, 1972.

\bibitem{BBM03}
F.~Borrelli, A.~Bemporad, and M.~Morari, ``Geometric algorithm for
  multiparametric linear programming,'' \emph{Journal of optimization theory
  and applications}, vol. 118, pp. 515--540, 2003.

\bibitem{BMDP02a}
A.~Bemporad, M.~Morari, V.~Dua, and E.~Pistikopoulos, ``The explicit linear
  quadratic regulator for constrained systems,'' \emph{Automatica}, vol.~38,
  no.~1, pp. 3--20, 2002.

\bibitem{TJB03}
P.~T{\o}ndel, T.~Johansen, and A.~Bemporad, ``An algorithm for multi-parametric
  quadratic programming and explicit {MPC} solutions,'' \emph{Automatica},
  vol.~39, no.~3, pp. 489--497, 2003.

\bibitem{GBN11}
A.~Gupta, S.~Bhartiya, and P.~Nataraj, ``A novel approach to multiparametric
  quadratic programming,'' \emph{Automatica}, vol.~47, no.~9, pp. 2112--2117,
  2011.

\bibitem{PS10}
P.~Patrinos and H.~Sarimveis, ``A new algorithm for solving convex parametric
  quadratic programs based on graphical derivatives of solution mappings,''
  \emph{Automatica}, vol.~46, no.~9, pp. 1405--1418, 2010.

\bibitem{patrinos2011convex}
------, ``Convex parametric piecewise quadratic optimization: Theory and
  algorithms,'' \emph{Automatica}, vol.~47, no.~8, pp. 1770--1777, 2011.

\bibitem{chen2022learning}
T.~Chen, X.~Chen, W.~Chen, H.~Heaton, J.~Liu, Z.~Wang, and W.~Yin, ``Learning
  to optimize: A primer and a benchmark,'' \emph{Journal of Machine Learning
  Research}, vol.~23, no. 189, pp. 1--59, 2022.

\bibitem{amos2023tutorial}
B.~Amos, ``Tutorial on amortized optimization,'' \emph{Foundations and Trends
  in Machine Learning}, vol.~16, no.~5, pp. 592--732, 2023.

\bibitem{sambharya2024learning}
R.~Sambharya, G.~Hall, B.~Amos, and B.~Stellato, ``Learning to warm-start
  fixed-point optimization algorithms,'' \emph{Journal of Machine Learning
  Research}, vol.~25, no. 166, pp. 1--46, 2024.

\bibitem{donti2021dc}
\BIBentryALTinterwordspacing
P.~L. Donti, D.~Rolnick, and J.~Z. Kolter, ``{DC}3: A learning method for
  optimization with hard constraints,'' in \emph{International Conference on
  Learning Representations}, 2021. [Online]. Available:
  \url{https://openreview.net/forum?id=V1ZHVxJ6dSS}
\BIBentrySTDinterwordspacing

\bibitem{grontas2026pinet}
\BIBentryALTinterwordspacing
P.~D. Grontas, A.~Terpin, E.~C. Balta, R.~D'Andrea, and J.~Lygeros, ``Pinet:
  Optimizing hard-constrained neural networks with orthogonal projection
  layers,'' in \emph{The Fourteenth International Conference on Learning
  Representations}, 2026. [Online]. Available:
  \url{https://openreview.net/forum?id=EJ680UQeZG}
\BIBentrySTDinterwordspacing

\bibitem{tordesillas2023rayen}
J.~Tordesillas, J.~P. How, and M.~Hutter, ``Rayen: Imposition of hard convex
  constraints on neural networks,'' \emph{arXiv preprint arXiv:2307.08336},
  2023.

\bibitem{ichnowski2021accelerating}
J.~Ichnowski, P.~Jain, B.~Stellato, G.~Banjac, M.~Luo, F.~Borrelli, J.~E.
  Gonzalez, I.~Stoica, and K.~Goldberg, ``Accelerating quadratic optimization
  with reinforcement learning,'' \emph{Advances in Neural Information
  Processing Systems}, vol.~34, pp. 21\,043--21\,055, 2021.

\bibitem{oshin2026deep}
\BIBentryALTinterwordspacing
A.~Oshin, R.~V. Ghosh, A.~D. Saravanos, and E.~Theodorou, ``Deep flex{QP}:
  Accelerated nonlinear programming via deep unfolding,'' in \emph{The
  Fourteenth International Conference on Learning Representations}, 2026.
  [Online]. Available: \url{https://openreview.net/forum?id=HL3TvE4Afm}
\BIBentrySTDinterwordspacing

\bibitem{langeMMOptimizationAlgorithms}
K.~Lange, ``{{MM Optimization Algorithms}} \textbar{} {{SIAM Publications
  Library}},'' https://epubs.siam.org/doi/book/10.1137/1.9781611974409.

\bibitem{bertsimas2022online}
D.~Bertsimas and B.~Stellato, ``Online mixed-integer optimization in
  milliseconds,'' \emph{INFORMS Journal on Computing}, vol.~34, no.~4, pp.
  2229--2248, 2022.

\bibitem{NEURIPS2020_6d0c9328}
\BIBentryALTinterwordspacing
K.~Wang, B.~Wilder, A.~Perrault, and M.~Tambe, ``Automatically learning compact
  quality-aware surrogates for optimization problems,'' in \emph{Advances in
  Neural Information Processing Systems}, H.~Larochelle, M.~Ranzato,
  R.~Hadsell, M.~Balcan, and H.~Lin, Eds., vol.~33.\hskip 1em plus 0.5em minus
  0.4em\relax Curran Associates, Inc., 2020, pp. 9586--9596. [Online].
  Available:
  \url{https://proceedings.neurips.cc/paper_files/paper/2020/file/6d0c932802f6953f70eb20931645fa40-Paper.pdf}
\BIBentrySTDinterwordspacing

\bibitem{wu2026koopman}
L.~Wu, W.~G.~Y. Tan, R.~D. Braatz, and J.~Drgo{\v{n}}a, ``Koopman-boxqp:
  Solving large-scale nmpc at khz rates,'' \emph{arXiv preprint
  arXiv:2602.18331}, 2026.

\bibitem{zhang2024latent}
Y.~Zhang, S.~Yang, T.~Ohtsuka, C.~Jones, and J.~Boedecker, ``Latent linear
  quadratic regulator for robotic control tasks,'' \emph{arXiv preprint
  arXiv:2407.11107}, 2024.

\bibitem{finn2017model}
C.~Finn, P.~Abbeel, and S.~Levine, ``Model-agnostic meta-learning for fast
  adaptation of deep networks,'' in \emph{International conference on machine
  learning}.\hskip 1em plus 0.5em minus 0.4em\relax PMLR, 2017, pp. 1126--1135.

\bibitem{rusu2018metalearning}
\BIBentryALTinterwordspacing
A.~A. Rusu, D.~Rao, J.~Sygnowski, O.~Vinyals, R.~Pascanu, S.~Osindero, and
  R.~Hadsell, ``Meta-learning with latent embedding optimization,'' in
  \emph{International Conference on Learning Representations}, 2019. [Online].
  Available: \url{https://openreview.net/forum?id=BJgklhAcK7}
\BIBentrySTDinterwordspacing

\bibitem{boyd2004convex}
S.~Boyd and L.~Vandenberghe, \emph{Convex optimization}.\hskip 1em plus 0.5em
  minus 0.4em\relax Cambridge university press, 2004.

\bibitem{durkan2019neural}
C.~Durkan, A.~Bekasov, I.~Murray, and G.~Papamakarios, ``Neural spline flows,''
  \emph{Advances in neural information processing systems}, vol.~32, 2019.

\bibitem{pmlr-v70-amos17b}
\BIBentryALTinterwordspacing
B.~Amos, L.~Xu, and J.~Z. Kolter, ``Input convex neural networks,'' in
  \emph{Proceedings of the 34th International Conference on Machine Learning},
  ser. Proceedings of Machine Learning Research, D.~Precup and Y.~W. Teh, Eds.,
  vol.~70.\hskip 1em plus 0.5em minus 0.4em\relax PMLR, 06--11 Aug 2017, pp.
  146--155. [Online]. Available:
  \url{https://proceedings.mlr.press/v70/amos17b.html}
\BIBentrySTDinterwordspacing

\bibitem{schallerLearningParametricConvex2025a}
M.~Schaller, A.~Bemporad, and S.~Boyd, ``Learning {{Parametric Convex
  Functions}},'' \emph{arXiv preprint}, no. arXiv:2506.04183, Jun. 2025.

\bibitem{molga2005test}
M.~Molga and C.~Smutnicki, ``Test functions for optimization needs,''
  \emph{Test functions for optimization needs}, vol. 101, no.~48, p.~32, 2005.

\bibitem{Bem25}
A.~Bemporad, ``An {L-BFGS-B} approach for linear and nonlinear system
  identification under $\ell_1$ and group-lasso regularization,'' \emph{IEEE
  Transactions on Automatic Control}, vol.~70, no.~7, pp. 4857--4864, 2025,
  code available at \url{https://github.com/bemporad/jax-sysid}.

\bibitem{kingma2014adam}
D.~P. Kingma and J.~Ba, ``Adam: A method for stochastic optimization,''
  \emph{arXiv preprint arXiv:1412.6980}, 2014.

\bibitem{byrd1995limited}
R.~H. Byrd, P.~Lu, J.~Nocedal, and C.~Zhu, ``A limited memory algorithm for
  bound constrained optimization,'' \emph{SIAM Journal on scientific
  computing}, vol.~16, no.~5, pp. 1190--1208, 1995.

\bibitem{qian2016hierarchical}
X.~Qian, A.~De~La~Fortelle, and F.~Moutarde, ``A hierarchical model predictive
  control framework for on-road formation control of autonomous vehicles,'' in
  \emph{2016 IEEE intelligent vehicles symposium (iv)}, 2016, pp. 376--381.

\bibitem{Andersson2018}
J.~A.~E. Andersson, J.~Gillis, G.~Horn, J.~B. Rawlings, and M.~Diehl,
  ``{CasADi} -- {A} software framework for nonlinear optimization and optimal
  control,'' \emph{Mathematical Programming Computation}, 2018.

\bibitem{Ferreau2014}
H.~Ferreau, C.~Kirches, A.~Potschka, H.~Bock, and M.~Diehl, ``{qpOASES}: A
  parametric active-set algorithm for quadratic programming,''
  \emph{Mathematical Programming Computation}, vol.~6, no.~4, pp. 327--363,
  2014.

\bibitem{agrawal2018rewriting}
A.~Agrawal, R.~Verschueren, S.~Diamond, and S.~Boyd, ``A rewriting system for
  convex optimization problems,'' \emph{Journal of Control and Decision},
  vol.~5, no.~1, pp. 42--60, 2018.

\bibitem{Clarabel_2024}
P.~J. Goulart and Y.~Chen, ``Clarabel: An interior-point solver for conic
  programs with quadratic objectives,'' \emph{arXiv preprint arXiv:2405.12762},
  2024.

\end{thebibliography}

\end{document}